# Delay-independent dual-rate PID controller for a packet-based Networked Control System

J. Alcaina[a], A. Cuenca[a,b,*], J. Salt[a], V. Casanova[a], R. Pizá[a]

[a] *Departamento de Ingeniería de Sistemas y Automática, Instituto Universitario de Automática e Informática Industrial, Universitat Politècnica de València, Spain*
[b] *Mechanical Systems Control Laboratory, Mechanical Engineering Department, University of California, Berkeley, United States*
joalac@upv.es, {acuenca, julian, vcasanov, rpiza}@isa.upv.es

## Abstract

In this paper, a novel delay-independent control structure for a networked control system (NCS) ~~which integrates packet-based control strategies with predictor-based and dual-rate control techniques, is proposed~~ is proposed, where packet-based control strategies with predictor-based and dual-rate control techniques are integrated. The control solution is able to cope with some networked communication problems such as time-varying delays, packet dropouts and packet disorder. In addition, the proposed approach enables to reduce ~~the~~ network load, and ~~the~~ usage of ~~the~~ connected devices, while maintaining a satisfactory control performance. As a delay-independent control solution, no network-induced delay measurement is needed for controller implementation. In addition, the control scheme is applicable to open-loop unstable plants. Control system stability is ensured in terms of linear matrix inequalities (LMIs). Simulation results show the main benefits of the control approach, which are experimentally validated by means of a Cartesian-robot-based test-bed platform.

* Author to whom correspondence should be addressed: acuenca@isa.upv.es.

| NOMENCLATURE | |
|---|---|
| $N, t, T, NT$ | Multiplicity, and basic, actuation and sensor periods |
| $\tau_k^{lr}, \tau_k^{rl}, \tau_k$ | Local-to-remote, remote-to-local, and round-trip time network-induced delay |
| $d_k^{lr}, d_k^{rl}$ | Local-to-remote and remote-to-local packet dropout occurrence |
| $y_k^{NT}, \hat{y}_k^{NT}, \bar{y}_k^{NT}$ | Actual, estimated, and actual or estimated, system output (at $NT$) |
| $r_k^{NT}$ | Reference signal (at $NT$) |
| $e_k^{NT}$ | Error signal (PI input) (at $NT$) |
| $x_k^{NT}, \hat{x}_k^{NT}, \bar{x}_k^{NT}$ | Actual, estimated, and actual or estimated, system state (at $NT$) |
| $v_k^{NT}, \hat{v}_k^{NT}, \bar{v}_k^{NT}$ | Actual, estimated, and actual or estimated, PI control action (at $NT$) |
| $\tilde{v}_k^{T}$ | Expanded PI control action (at $T$) |
| $v_k^{T}, \hat{v}_k^{T}, \bar{v}_k^{T}$ | Actual, estimated, and actual or estimated, PI action converted at $T$ |
| $u_k^{T}, \hat{u}_k^{T}, \bar{u}_k^{T}$ | Actual, estimated, and actual or estimated PD control action (at $T$) |
| $M$ | Upper bound for estimations |
| $z$ | $t$-unit operator |
| $\bar{z} = z^L, \quad L = \dfrac{T}{t}$ | $T$-unit operator |
| $z_N = \bar{z}^N = z^{N \cdot L}$ | $NT$-unit operator |

# 1. - Introduction

Networked Control Systems (NCSs) [33, 36] is a prolific control area, addressing control scenarios where different devices share a common communication link. There are several advantages associated with NCS such as cost reduction, ~~flexibility~~ and ease of installation and maintenance, but also drawbacks like the possible occurrence of time-varying delays [12, 26, 28-32, 37, 40], packet dropouts [9, 13, 18, 26, 30], packet disorder [9, 19, 20, 26, 30], and network bandwidth constraints [5, 11, 16]. One interesting aim in NCS is to reduce the number of transmissions through the network, which can result in some advantages such as increase in ~~the~~ network bandwidth, and enlargement of ~~the~~ battery usage of the different wireless devices connected to the NCS. But the reduction in transmissions should come along with preservation of stability and performance properties. In order to ensure this achievement, different control solutions have been proposed in literature, for instance: packet-based control [38, 39, 41], which enables to decrease the communication rate by simultaneously sending a set of data in each transmission; event-based control [4, 15, 32, 34, 35, 42], where the transmission is only carried out if control or output variables satisfy a certain event condition; multi-rate control [6, 8, 24], where a slower sensing rate in comparison to a faster actuation can be assumed; and predictor-based control [10, 14, 27, 35], which exploits model-based predictions to address scarce data and compensate for network-induced delays. The present work proposes a control structure for an NCS, which integrates packet-based control strategies with predictor-based and dual-rate control techniques. The systematic combination of these control solutions enables not only to reduce the number of transmissions through the network but also to face some networked communication problems such as time-varying delays, packet dropouts, and packet disorder. In addition, the control approach is able to keep a satisfactory control performance, which is defined by means of the nominal (no-delay, no-dropout) dual-rate response.

Dual-rate control techniques provide twofold benefits: to avoid packet disorder, and to reduce the number of transmissions through the network. The sensing period may be chosen to be larger than the largest round-trip time delay found in a statistical distribution for the network-induced delay, which is assumed to be known [2, 23]. In this way, no packet disorder is produced, and the use of the network and devices can be reduced. Despite the fact of choosing a slow sensing period, a satisfactory control performance may be achieved by selecting a faster actuation period [24]. In the present work, due to the wide knowledge of PID controllers in industrial and academic environments, a dual-rate PID control structure is taken into account. The controller is split into two parts: a slow-rate PI controller and a fast-rate PD controller. The integral action is applied at slow rate, because it usually operates at this frequency zone. As data travel through the network at slow rate, the PI controller is located at the remote side, with no direct connection to the plant. The derivative actions, which are associated with faster dynamics, are applied at fast rate to reach the satisfactory control performance, and hence, the PD controller is located at the local side, directly connected to the plant. Basic design procedures of dual-rate PID controllers can be looked up in [24].

Combining packet-based control strategies with predictor-based control techniques enables to deal with packet dropouts and time-varying delays. Packet dropouts may occur because of the utilization of user datagram protocol (UDP) as the communication protocol in this work [2, 23]. Predictor-based control techniques can be used at the remote side in order to deal with up to $M$ possible consecutive packet dropouts, being $M$ an upper-bound which may be established from off-line experiences ~~on~~ in the NCS. The *M* future, estimated PI control actions will be sent to the local side by implementing a packet-based control strategy in order to compute the next PD control actions. This computation can be carried out following a delay-free control algorithm, which is a central aspect of this work. At the local side, when no packet arrives, the PD controller is able to compute an estimated control signal from the PI control actions received in the previous transmission, and then, it can apply the signal following a uniform actuation pattern, that is, not being influenced by the time-varying network-induced delay. When the packet arrives after the delay, the PD controller is able to compute the actual control actions irrespective of the delay, and replace the estimated control signal with the actual one. Assuming neither model uncertainties nor disturbances, the difference between actual and estimated control signals should be negligible. Note that, inside the current sensing period, the control signal estimated for this period is injected from the beginning of the period to the moment in which the new packet is received. This is an essential difference compared with [7], where the last control action of the previous sensing period is held along the current sensing period until the new packet arrives, being additionally required a gain-scheduling approach to retune the controller according to the delay. To the best of the authors' knowledge, the working mode of the delay-independent control algorithm proposed in this work is novel in this kind of frameworks. Since the present work is able to consider unstable plants to be controlled, the state prediction solution will include a state resetting procedure [7, 21, 22].

In summary, the proposed delay-independent dual-rate PID controller may be defined as a PID controller, which is able to generate a fast-rate control signal from a slow-rate process output measurement, calculating the fast-rate control signal by following a delay-free control algorithm, which means that the network-induced delay does not affect the computation. Therefore, the motivation of introducing this kind of controller is twofold: i) for implementation purposes, an obvious advantage derived from the delay-independent control solution is that the round-trip time delay is not required to be measured. This feature makes the solution applicable to a wide range of NCSs where the time delay is difficult to measure; and ii) no additional control techniques (e.g., gain-scheduling control [7, 25], optimal control [19], fuzzy control [31, 32], H∞ filtering techniques [30], etc.) may be required to compensate for the delay, which reduces the complexity of the control solution.

Although the proposed controller is delay-independent, the plant, due to the time-varying network delay, is subjected to some variations in the instants where the input commands are presented to the plant. Then, the plant model can vary from sensor period to sensor period. If that were not the case, the eigenvalues of the NCS closed-loop model would determine stability and performance, but eigenvalues are of limited usefulness in time-varying contexts, such as the networked one in this paper. This leads to represent the NCS closed-loop model as a Linear Time Varying (LTV) system, which depends on the time-varying delay, and whose stability can be ensured via Linear Matrix Inequalities (LMIs).

Summarizing, the main contribution of this paper is the development of a new and comprehensive approach, where dual-rate and predictor-based control techniques and packet-based control strategies are systematically brought together in an NCS in order to face some communication drawbacks and reduce ~~the~~ resource utilization, while maintaining a satisfactory control performance. The central feature of the proposed control solution is its delay-free control signal computation. This is a distinguished improvement compared to other control solutions, where the network-induced delay must be compensated for, and hence, additional control techniques may be required.

The paper is structured in the following sections. In section 2, the problem scenario is formally introduced. In section 3, the control techniques used at the remote and local sides are presented. Control system stability is enunciated in terms of LMIs in section 4. Simulation results in section 5 illustrate the benefits of the proposed control strategy in an unstable open-loop plant; concretely, by controlling the position output of a Cartesian robot. In Section 6, the previous results are experimentally validated in a test-bed platform based on the Cartesian robot, and using UDP as the network protocol. Finally, conclusions close this contribution.

## 2. - Problem description

The proposed NCS is depicted in Figure 1, where the network is placed between the remote and local sides, and it introduces some communication problems such as time-varying delays, packet dropouts, and packet disorder. In the next subsections, these problems are formally described, and the control structure is presented in detail.

*2.1. Time-varying delays, packet dropouts, and packet disorder*

The round-trip time delay for the packet sampled at instant $kNT$ (where $k \in \mathbb{N}$, $T$ is the actuation period, $NT$ is the sensor period, and $N \in \mathbb{Z}^+$ is a parameter known as multiplicity in a dual-rate control framework [24]) is defined as

$$\tau_k = \tau_k^{lr} + \tau_k^{rl} + \tau_k^{c}, \tag{2.1}$$

being $\tau_k^{lr}$ the local-to-remote network-induced delay, $\tau_k^{rl}$ the remote-to-local delay, and $\tau_k^{c}$ a negligible computation time delay. In order to avoid packet disorder, in this work it is strictly necessary to know the maximum round-trip time delay $\tau_{\max}$ such that $\tau_{\max} < NT$. As a common timer is supposed to be shared by the local devices in such a way that all of them are perfectly synchronized, $\tau_k$ can be measured subtracting packet sending and receiving times, not requiring time-stamping techniques. Therefore, from off-line experiences ~~on~~ in the NCS, the statistical distribution of the round-trip time delay can be obtained, and hence, $\tau_{max}$. Since in this work, an IP network which uses UDP as the transport layer protocol is taken into account, the distribution will be a constant plus a Gamma distributed random variable, whose shape and scale parameters change with load and network segment [23]. Usually, this distribution is approximated as a generalized exponential distribution [25], whose probability density function can take this form

$$P[\tau_k] = \begin{cases} \frac{1}{\phi} e^{\frac{-(\tau_k-\eta)}{\phi}}, & \tau_k \geq \eta \\ 0, & \tau_k < \eta \end{cases}, \quad (2.2)$$

being the expected value of the delay $E[\tau_k] = \phi + \eta$, and its variance $V[\tau_k] = \phi^2$. A feasible choice of $\eta$ is the median of the delay. From $\eta$ and an experimental value of $E[\tau_k]$ (or the mean), $\phi$ can be easily approximated.

As well-known, when using the UDP transmission model, packet dropouts may appear. This phenomenon is essentially random [2], and hence, it can be modeled as a Bernoulli distribution [33]. The variable $d_k^{lr}$ indicates the possible loss of the packet sent from the local side to the remote side at instant *kNT*. Similarly, $d_k^{rl}$ is defined for the opposite network link. In this work, both variables are considered as a Bernoulli process with probability of dropout:

$$\begin{aligned} p^{lr} &= \Pr[d_k^{lr} = 0] \in [0,1) \\ p^{rl} &= \Pr[d_k^{rl} = 0] \in [0,1) \end{aligned}, \quad (2.3)$$

In some real scenarios, $p^{lr}, p^{rl}$ could be considered as the same value $p = p^{lr} = p^{rl}$.

Finally, let us define *M* as the upper bound of consecutive packet dropouts. From a significant number of off-line experiences in the considered NCS, where a certain probability of dropout for the Bernoulli process in (2.3) can be chosen, one can observe the resulting number of consecutive packet dropouts which is produced in each experience. The largest number of consecutive packet dropouts for the whole set of experiences can be determined as the worst case. Although this number might appear with low probability, in order to consider a conservative decision it may be assigned to *M*.

*2.2. Control structure*

Next, the different devices included in Figure 1 for the NCS are presented:

- the process to be controlled is a Cartesian robot. More information about this process can be found in Sections 5 and 6.

- the sensor works at period *NT* to sample the process output $y_k^{NT}$, which, in this case, is the position of the robot. Sensing at this slow rate enables to reduce the number of transmissions, reducing network load and device usage.

- the slow-rate PI controller generates a PI control action $v_k^{NT}$ from the reference $r_k^{NT}$ and the output $y_k^{NT}$, as long as the output arrives to the remote side (when $d_k^{lr}=1$) after $\tau_k^{lr}$. Otherwise (when $d_k^{lr}=0$), a previously estimated PI control action $\hat{v}_k^{NT}$ will be used. Let us consider a maximum waiting time $\tau_{\max}^{lr}$ to detect a packet dropout in this device. If $\tau_{\max}^{lr}$ expires and the packet does not arrive, it will be considered as a dropout. More information about the definition and operation mode of the slow-rate PI controller can be found in subsection 3.2.

- the prediction stage computes an array of *M* estimated, future PI control actions $\left[\hat{v}_{k+1}^{NT},\hat{v}_{k+2}^{NT},\ldots,\hat{v}_{k+M}^{NT}\right]$ taking into account: i) the array of the actual and future references $\left[r_k^{NT},r_{k+1}^{NT},r_{k+2}^{NT},\ldots r_{k+M}^{NT}\right]$, ii) the actual PI control action $v_k^{NT}$ or its estimation $\hat{v}_k^{NT}$, iii) the actual process output $y_k^{NT}$ or its estimation $\hat{y}_k^{NT}$, and iv) the actual state $x_k^{NT}$ or its estimation $\hat{x}_k^{NT}$. For the sake of simplicity and brevity, both cases (the actual and estimated ones) will be contained under the notation $\bar{v}_k^{NT},\bar{y}_k^{NT},\bar{x}_k^{NT}$ in the sequel. The main goals of the prediction stage are: i) to face packet dropouts for both network links; and ii) to provide estimated PI control actions in order to be suitably used at the local side to compute the PD control signal via a delay-independent control algorithm. For more information about how the prediction stage works, see subsection 3.5.

- the packet generator implements a packet-based control strategy, which creates the packet to be sent to the local side, containing the set of estimated PI control actions $\left[\bar{v}_k^{NT}, \hat{v}_{k+1}^{NT}, \hat{v}_{k+2}^{NT}, \ldots, \hat{v}_{k+M}^{NT}\right]$.

- the actuator may include a rate converter and a fast-rate PD controller: firstly, the rate converter converts the slow-rate PI control signal into a fast-rate one to be used by the fast-rate PD controller as an input (more details in subsection 3.3). Secondly, the controller generates the fast-rate PD control signal in order to achieve the desired performance, which is defined by the nominal (no-delay, no-dropout) dual-rate control strategy. According to $d_k^{rl}$, the actuator applies the PD control actions using a different actuation pattern:

a) if dropout occurs ($d_k^{rl} = 0$): the actuator injects the control actions following a uniform pattern at time instants $\{0, T, 2T, \ldots, (N-1)T\}$ inside the current sensor period $NT$ (see Figure 2). In this case, the PD control signal is computed from an estimated PI control value received in a previous successful delivery. Then, an estimated PD control signal $\left[\hat{u}_k^T, \hat{u}_{k+1}^T, \ldots, \hat{u}_{k+N-1}^T\right]$ is generated.

b) if no dropout occurs ($d_k^{rl} = 1$): the actuator injects the control actions following a non-uniform pattern. For a particular pattern where $\tau_k < T$, the actuation time instants would be $\{0, \tau_k, T, 2T, \ldots, (N-1)T\}$ inside the current sensor period $NT$ (see Figure 3). In this case, the PD control signal would take this form $\left[\hat{u}_k^T, \bar{u}_k^T, \bar{u}_{k+1}^T, \ldots, \bar{u}_{k+N-1}^T\right]$, which would mean the injection of an estimated PD control action $\hat{u}_k^T$ at the beginning of the current sensor period, and $N$ actual or estimated PD control actions $\bar{u}_{k+i}^T, i = [0..N-1]$ at the rest of the time instants after $\tau_k$. This subset of $N$ control actions will be composed of actual values, when an actual PI control signal $v_k^{NT}$ is used to compute them. Otherwise, when an estimated PI control signal $\hat{v}_k^{NT}$ is used, the subset will be composed of estimated values. In any case, assuming neither model uncertainties nor disturbances, $\hat{u}_k^T$ should be very similar to $\bar{u}_k^T$.

Under both actuation patterns, the PD control actions can be computed irrespective of the delay, resulting in a delay-independent control solution. This important feature distinguishes the solution, since no delay measurement is needed to implement the controller. More information about the fast-rate PD controller can be found in subsection 3.4.

## 3. – Control design.

The control proposal will be formulated in the next subsections. Firstly, some preliminaries, which are needed for the design step, will be presented. Secondly, the dual-rate controller will be stated. As commented in section 2, the controller is composed of a slow-rate PI controller, a fast-rate PD controller, and between them, a rate converter. Subsections 3.2 to 3.4 will detail the design aspects for each part of the dual-rate controller, differentiating between the dropout case and the no-dropout case. Finally, in subsection 3.5, the prediction stage will be enunciated.

*3.1. Preliminaries*

Let us define the transfer function of the continuous plant to be controlled as $G_p(s)$. By using the Z-transform at different periods plus a zero-order hold device $H(s)$, different discrete-time versions for $G_p(s)$ can be considered. Therefore, denoting $\bar{z} = z^L$ and $z_N = \bar{z}^N = z^{N \cdot L}$:

$$G^{NT}(z_N) \triangleq Z_{NT}\left[H_{NT}G_p(s)\right] = \frac{Y^{NT}(z_N)}{U^{NT}(z_N)};\ G^T(\bar{z}) \triangleq Z_T\left[H_T G_p(s)\right] = \frac{Y^T(\bar{z})}{U^T(\bar{z})}$$

$$G^t(z) \triangleq Z_t\left[H_t G_p(s)\right] = \frac{Y^t(z)}{U^t(z)}, \quad t < T: \ t \cdot L = T, \quad L \in \mathbb{Z}^+,$$

(3.1)

In addition, the consequent state-space representations for each case (where matrices have suitable dimensions) can be enunciated as

$$\begin{cases} x_{k+1}^{NT} = A^{NT} x_k^{NT} + B^{NT} u_k^{NT} \\ y_k^{NT} = C^{NT} x_k^{NT} \end{cases} ; \begin{cases} x_{k+1}^{T} = A^{T} x_k^{T} + B^{T} u_k^{T} \\ y_k^{T} = C^{T} x_k^{T} \end{cases} ; \begin{cases} x_{k+1}^{t} = A^{t} x_k^{t} + B^{t} u_k^{t} \\ y_k^{t} = C^{t} x_k^{t} \end{cases}, \quad (3.2a)$$

Finally, let us consider a continuous PID controller, which is designed according to classical methods in order to achieve some specifications for the process to be controlled. The controller gains used for the next design steps are given by this PID configuration:

$$G_{PID}(s) = K_p \left(1 + T_d s + \frac{1}{T_i s}\right), \quad (3.2b)$$

*3.2. Slow-rate PI controller*

a) No-dropout case ($d_k^{lr} = 1$): The PI controller working at period *NT* is enunciated as

$$G_{PI}^{NT}(z_N) = K_{PI} \frac{z_N - \left(1 - \dfrac{NT}{T_i}\right)}{z_N - 1} = \frac{V^{NT}(z_N)}{E^{NT}(z_N)}, \quad (3.3)$$

being $V^{NT}(z_N)$ the PI control signal, $E^{NT}(z_N)$ the error signal, and $K_{PI}, T_i$ the gains of the PI controller. Let us consider $K_{PI} = 1$. The PI control signal is obtained as

$$V^{NT}(z_N) = G_{PI}^{NT}(z_N) E^{NT}(z_N) = G_{PI}^{NT}(z_N) \left(R^{NT}(z_N) - Y^{NT}(z_N)\right), \quad (3.4)$$

and, from (3.3), the difference equation for the PI controller with $K_{PI} = 1$ will be

$$v_k^{NT} = v_{k-1}^{NT} + e_k^{NT} - \left(1 - \frac{NT}{T_i}\right) e_{k-1}^{NT} = v_{k-1}^{NT} + \left(r_k^{NT} - y_k^{NT}\right) - \left(1 - \frac{NT}{T_i}\right) \left(r_{k-1}^{NT} - y_{k-1}^{NT}\right), \quad (3.5)$$

b) Dropout case ($d_k^{lr} = 0$): In this case, instead of using the actual PI control signal in (3.4), the estimated one $\hat{V}^{NT}(z_N)$ must be used. This signal is previously generated at the prediction stage according to subsection 3.5.

*3.3. Rate converter*

As it is well-known [24], a rate converter $\left[H^{NT}\right]^T(\bar{z})$ is required between slow and fast rate controllers. Its goal is to convert the slow-rate PI control signal into a fast-rate one to be used by the fast-rate PD controller as an input. This operation, with low computational complexity, can be carried out at the local side. Two cases are considered:

a) No-dropout case ($d_k^{lr} = 1$): The rate converter considers the actual slow-rate PI control signal $V^{NT}(z_N)$ to obtain the held fast-rate one $V^T(\bar{z})$. As, in this work, step references are considered, the rate converter becomes a digital zero-order hold:

$$\left[H^{NT}\right]^T(\bar{z}) = \frac{V^T(\bar{z})}{\left[V^{NT}(z_N)\right]^T} = \frac{1-\bar{z}^{-N}}{1-\bar{z}^{-1}} \rightarrow V^T(\bar{z}) = \left[H^{NT}\right]^T(\bar{z}) \cdot \left[V^{NT}(z_N)\right]^T, \quad (3.6a)$$

Note that $V^{NT}(z_N)$ is required to be used in an expanded way $\left[V^{NT}(z_N)\right]^T$, that is,

$$\left[V^{NT}(z_N)\right]^T \triangleq \tilde{V}^T(\bar{z}) \triangleq \sum_{k=0}^{\infty} \tilde{v}_k^T \bar{z}^{-k} : \begin{cases} \tilde{v}_k^T = v_k^T, \forall k = \lambda N \\ \tilde{v}_k^T = 0, \quad \forall k \neq \lambda N \end{cases}, \lambda \in Z^+, \quad (3.6b)$$

More information can be found in [24].

b) Dropout case ($d_k^{lr} = 0$): Now, the rate converter considers the estimated PI control signal $\hat{V}^{NT}(z_N)$

$$\hat{V}^T(\bar{z}) = \left[H^{NT}\right]^T(\bar{z}) \cdot \left[\hat{V}^{NT}(z_N)\right]^T, \quad (3.7)$$

As said in section 2, for the sake of simplicity and brevity, both cases ((3.6a) and (3.7)) will be contained under the notation $\bar{V}^T(\bar{z})$ from now on.

*3.4. Fast-rate PD controller*

a) No-dropout case ($d_k^{rl} = 1$): Let us define $K_{PD} = K_p$. Then, the controller is stated as

$$G_{PD}^T(\bar{z}) = K_{PD} \frac{\bar{z}\left(1+\frac{T_d}{T}\right) - \frac{T_d}{T}}{\bar{z}} = \frac{\bar{U}^T(\bar{z})}{\bar{V}^T(\bar{z})} \rightarrow \bar{U}^T(\bar{z}) = G_{PD}^T(\bar{z})\bar{V}^T(\bar{z}), \quad (3.8)$$

and its difference equation is

$$\bar{u}_k^T = K_{PD}\left(1+\frac{T_d}{T}\right)\bar{v}_k^T - K_{PD}\left(\frac{T_d}{T}\right)\bar{v}_{k-1}^T, \quad (3.9)$$

Observe that the notation $\bar{U}^T(\bar{z})$ represents both $U^T(\bar{z})$ and $\hat{U}^T(\bar{z})$, where $U^T(\bar{z})$ is the actual PD control signal obtained from the actual PI control signal $V^T(\bar{z})$ (in (3.6a)), and $\hat{U}^T(\bar{z})$ is the estimated PD control signal obtained from the estimated PI control signal $\hat{V}^T(\bar{z})$ (in (3.7)). Iterating the difference equation deduced from (3.9) $N$ times, the $N$ PD control actions are generated and applied after $\tau_k$, which is the moment when $\bar{V}^T(\bar{z})$ is available. As commented in Section 2, due to the delay, these actions will be applied following a non-uniform pattern. As there are different patterns depending on $\tau_k$, for the sake of clarity let us consider the particular case where $\tau_k < T$ (as in Section 2). Then, a basic period $t$ is required to adapt the non-uniformity to the delay in such a way that the actuation pattern inside the sensor period $NT$ will take this form (where $l=0..LN\text{-}1$):

$$\begin{cases} \hat{u}_k^T, & lt = 0..\tau_k \\ \bar{u}_k^T, & lt = \tau_k..T \\ \bar{u}_{k+1}^T, & lt = T..2T \\ \vdots \\ \bar{u}_{k+N-1}^T, & lt = (N-1)T..NT \end{cases}, \quad (3.10)$$

Note that the control action $\hat{u}_k^T$, which is injected at the beginning of the sensor period $NT$, is an estimated control action. Depending on the occurrence, or not, of a dropout in the local-to-remote link, the rest of control actions in (3.10), $\bar{u}_{k+i}^T, i=0..N-1$, will be estimated actions $\hat{u}_{k+i}^T$ or actual actions $u_{k+i}^T$, respectively. In the first case, $\bar{u}_k^T = \hat{u}_k^T$, and hence (3.10) is equivalent to (3.12). In the second case, $\bar{u}_k^T = u_k^T$, and assuming an accurate prediction, the difference between $\hat{u}_k^T$ and $u_k^T$ will be negligible.

b) Dropout case ($d_k^{rl} = 0$): The estimated PI control signal $\hat{V}^T(\bar{z})$ (to be defined in the last step in subsection 3.5) is now required. This control signal is available at the local side, since it was received in a previous successful communication.

The estimated PD control signal takes this form

$$\hat{U}^T(\bar{z}) = G_{PD}^T(\bar{z})\hat{V}^T(\bar{z}), \qquad (3.11)$$

From (3.9), but considering now estimated signals, the set of *N* control actions are computed and applied according to the next uniform actuation pattern inside the sensor period *NT*:

$$\begin{cases} \hat{u}_k^T, & kT = 0..T \\ \hat{u}_{k+1}^T, & kT = T..2T \\ \vdots \\ \hat{u}_{k+N-1}^T, & kT = (N-1)T..NT \end{cases}, \qquad (3.12)$$

## 3.5. Prediction stage

The prediction algorithm is executed *M* times (*M* was defined in section 2 as the upper bound of consecutive packet dropouts) in order to generate the packet that includes the *M* future, estimated PI control actions $\left[\hat{v}_{k+1}^{NT}, \hat{v}_{k+2}^{NT}, \ldots, \hat{v}_{k+M}^{NT}\right]$. This packet is computed for every sensor period at the remote side, and it is sent to the local side in order to be used when subsequent dropouts occur through the remote-to-local communication link. Considering a for-loop, where *i*=1..*M*, the statements of the prediction algorithm included in the loop are based on the next steps:

a. Resetting of the initial state: If the current state sensed at period *NT*, $x_k^{NT}$, is available at the remote side, that is, no dropout occurs when being sent via the local-to-remote network link, a resetting of the initial condition for the state at period *t* and at period *T* can be carried out. This operation can be executed when *i*=1, and it is required to deal with possible unstable plants [7, 21, 22]. For the rest of iterations of the algorithm (*i*=2..*M*), or if the current state was dropped when *i*=1, the updating is computed from the estimated state $\hat{x}_{k+i-1}^{NT}$, which will be defined in step 3. As in section 2, to contemplate every situation, let us define a generic (actual or estimated) state $\bar{x}_k^{NT}$. Therefore, the resetting carried out in each iteration is

$$\begin{cases} i=1: & \hat{x}_k^T \leftarrow \bar{x}_k^{NT}; \quad \hat{x}_k^t \leftarrow \bar{x}_k^{NT} \\ i>1: & \hat{x}_{k+(i-1)N}^T \leftarrow \hat{x}_{k+i-1}^{NT}; \quad \hat{x}_{k+(i-1)LN}^t \leftarrow \hat{x}_{k+i-1}^{NT} \end{cases}, \quad (3.13)$$

b. Estimation of the *N* PD control actions either from the estimated PI control signal $\hat{V}^T(\bar{z})$ (this case can occur when $i \geq 1$) or from the actual one $V^T(\bar{z})$ (this case can occur only when *i*=1). Both cases assume the rate conversion previously carried out in (3.7) or (3.6a), respectively.

Similarly to (3.9), the estimated control signal is computed by iterating the next difference equation for *j*=0..*N*-1. Each iteration *i* for the prediction algorithm is calculated as follows

$$\begin{cases} i=1: & \hat{u}_{k+j}^T = K_{PD}\left(1+\dfrac{T_d}{T}\right)\bar{v}_{k+j}^T - K_{PD}\left(\dfrac{T_d}{T}\right)\bar{v}_{k-1+j}^T \\ i>1: & \hat{u}_{k+j+(i-1)N}^T = K_{PD}\left(1+\dfrac{T_d}{T}\right)\hat{v}_{k+j+(i-1)N}^T - K_{PD}\left(\dfrac{T_d}{T}\right)\hat{v}_{k-1+j+(i-1)N}^T \end{cases}, \quad (3.14)$$

c. Estimation of the next state and output at period *NT* from the estimated PD control actions. Now, as in (3.12), a uniform pattern is used. Then, for each iteration of the prediction algorithm *i*=1..*M*, the next state-space representation at period *T* is computed for *j*=0..*N*-1:

$$\begin{cases} \hat{x}_{k+1+j+(i-1)N}^T = A^T \hat{x}_{k+j+(i-1)N}^T + B^T \hat{u}_{k+j+(i-1)N}^T \\ \hat{y}_{k+1+j+(i-1)N}^T = C^T \hat{x}_{k+1+j+(i-1)N}^T \end{cases}, \quad (3.15)$$

As a result of iterating (3.15) for all *i*, the *M* estimated states and outputs at period *NT*, $\hat{x}_{k+i}^{NT}, \hat{y}_{k+i}^{NT}$, are calculated.

d. Estimation of the PI control signal $\hat{V}^{NT}(z_N)$ from the estimated output signal $\hat{Y}^{NT}(z_N)$. Note that, particularly for the first iteration of the prediction algorithm (*i*=1), the actual output $y_k^{NT}$ can be used if it is available at the remote side, that is, if no dropout occurs when being sent through the local-to-remote network link ($d_k^{lr}=1$). Then, the actual PI control action $v_k^{NT}$, which is generated by the output $y_k^{NT}$ (remember (3.5)), can also be used. In this way, similarly to step 1, a resetting of the initial condition for the PI controller ($v_k^{NT}$) is carried out in order to compute the next estimated PI control action $\hat{v}_{k+1}^{NT}$. This operation is useful due to the unstable open-loop nature of the PI controller [7, 21, 22]. As usual, in order to contemplate every situation in the prediction algorithm, let us define a generic (actual or estimated) output $\bar{y}_k^{NT}$, and a generic (actual or estimated) control action $\bar{v}_k^{NT}$. Therefore, similarly to (3.5), the different iterations $i$ of the prediction algorithm take the form

$$\begin{cases} i = 1: \quad \hat{v}_{k+1}^{NT} = \bar{v}_k^{NT} + \left(r_{k+1}^{NT} - \hat{y}_{k+1}^{NT}\right) - \left(1 - \frac{NT}{T_i}\right)\left(r_k^{NT} - \bar{y}_k^{NT}\right) \\ i > 1: \quad \hat{v}_{k+i}^{NT} = \hat{v}_{k+i-1}^{NT} + \left(r_{k+i}^{NT} - \hat{y}_{k+i}^{NT}\right) - \left(1 - \frac{NT}{T_i}\right)\left(r_{k+i-1}^{NT} - \hat{y}_{k+i-1}^{NT}\right) \end{cases}, \qquad (3.16)$$

## 4. Stability analysis

### 4.1 Closed-loop model via lifting

Let us consider a continuous linear time-invariant plant, which admits a state-space realization $\Sigma = (\bar{A}, \bar{B}, C, D)$, with suitable dimensions. Being $\xi$ an arbitrary real number, one can denote $B(\xi) = \int_0^{\xi} e^{\bar{A}\gamma} \bar{B} d\gamma$ if $\xi > 0$, or $B(\xi) = 0$ if $\xi \leq 0$. The discrete time sampled-data model at period $T$ of the previous plant was presented in (3.2a), being $A^T = e^{\bar{A}T}; \quad B^T = B(T); \quad C^T = C$.

In order to reflect the dual-rate sampling with either uniform actuation (when $d_k^{rl} = 0$) or non-uniform actuation (when $d_k^{rl} = 1$), the control system can be modeled via state-space representation adopting the so-called lifting methodology [17]. Then, let us represent the process as:

$$\begin{cases} x_{k+1}^{NT} = A_p x_k^{NT} + B_p \tilde{U}_k^T \\ y_k^{NT} = C_p x_k^{NT} \end{cases}, \qquad (4.1)$$

where

- $A_p = e^{\bar{A}NT}$; $B_p = \begin{bmatrix} B_0^* & B_1^* & \cdots & B_{\bar{N}}^* \end{bmatrix}$; $C_p = C$; where $\bar{N} = N$ for the non-uniform case, and $\bar{N} = N-1$ for the uniform case. In addition, $B_i^* = B(\lambda_{i+1} - \lambda_i)e^{\bar{A}(NT-\lambda_{i+1})}$, where $\lambda_i \ (i=0,1,\ldots,\bar{N})$ are the actuation time instants.

- $\tilde{U}_k^T$ is a generic array of control signals to indistinctly represent either the uniform actuation as in (3.12) via $\hat{U}_k^T$, or the non-uniform actuation as in (3.10) via $\hat{\bar{U}}_k^T$, where the particular case $\tau_k < T$ is considered.

The predictor stage can be modelled as follows:

$$\hat{x}_{k+1}^{NT} = \begin{cases} A_p \hat{x}_k^{NT} + B_p \tilde{U}_k^T, & \text{if } d_k^{lr} = 0 \\ A_p x_k^{NT} + B_p \tilde{U}_k^T, & \text{if } d_k^{lr} = 1 \end{cases}, \quad (4.2)$$

For the sake of simplicity and brevity, let us consider the following assumptions to define the state-space representation of the delay-independent controller:

1) the setpoints are constant, and hence we can assume them to be zero without loss of generality. Then, $e_k^{NT} = -y_k^{NT}$.

2) the ~~behaviour~~ behavior of the control system when $d_k^{lr} = 0$ is similar to the one when $d_k^{rl} = 0$, since in both cases the control signal $\tilde{U}_k^T$ is computed from the estimated output $\hat{y}_k^{NT}$.

3) the dual-rate controller can be defined as a cascade-connected system.

~~4) the network delay fulfills $\tau_k < T$ as in Section 5 in order to check stability for the control system in that section.~~

After manipulating the difference equations (3.5) and (3.9), and assuming that ~~the~~ integral actions are generated at slow rate, and ~~the~~ derivative actions are computed at fast rate (as commented in Section 1), the dual-rate controller can take this lifted expression

$$\begin{cases} \varphi_{k+1}^{NT} = A_c \varphi_k^{NT} - B_c \hat{y}_{k+1}^{NT} + \bar{B}_c \hat{y}_k^{NT} \\ \hat{U}_k^T = C_c \varphi_k^{NT} \end{cases}, \text{if } d_k^{lr} = 0 \text{ or } d_k^{rl} = 0$$
$$\begin{cases} \varphi_{k+1}^{NT} = A_c \varphi_k^{NT} - B_c y_{k+1}^{NT} + \bar{B}_c y_k^{NT} \\ \hat{\bar{U}}_k^T = \chi U_k^T + (\chi_I - \chi)\hat{U}_k^T = \bar{C}_c \varphi_k^{NT} \end{cases}, \text{if } d_k^{lr} = d_k^{rl} = 1$$
(4.3)

where, using $(\cdot)^T$ as transpose function, the controller state is $\varphi_k^{NT} = \begin{pmatrix} v_k^{NT} & \mu_k^{NT} \end{pmatrix}^T$, being $\mu_k^{NT} = v_{k-1}^{NT}$ and

$$A_c = \begin{pmatrix} 1 & 0 \\ 1 & 0 \end{pmatrix}; \quad B_c = \begin{pmatrix} 1 \\ 0 \end{pmatrix}; \quad \bar{B}_c = \begin{pmatrix} 1 - \dfrac{NT}{T_i} \\ 0 \end{pmatrix};$$

$$C_c = \begin{pmatrix} K_{PD}\left(1 + \dfrac{T_d}{T}\right) & -K_{PD}\dfrac{T_d}{T} \\ K_{PD} & 0 \\ \vdots & \vdots \\ K_{PD} & 0 \end{pmatrix}_{N \times 2},$$
(4.4)

In addition, depending on the delay $\tau_k$, different configurations for $\chi$, $\chi_I$ and $\bar{C}_c$ may be considered. Therefore, when $\tau_k < T$

$$\chi = \begin{pmatrix} 0_{1 \times N} \\ I_{N \times N} \end{pmatrix}; \quad \chi_I = \begin{pmatrix} I_{1 \times N} \\ I_{N \times N} \end{pmatrix};$$

$$\bar{C}_c = \begin{pmatrix} K_{PD}\left(1 + \dfrac{T_d}{T}\right) & -K_{PD}\dfrac{T_d}{T} \\ K_{PD}\left(1 + \dfrac{T_d}{T}\right) & -K_{PD}\dfrac{T_d}{T} \\ K_{PD} & 0 \\ \vdots & \vdots \\ K_{PD} & 0 \end{pmatrix}_{(N+1) \times 2},$$
(4.5a)

When $\tau_k > dT$, $d = 1..N-1$

$$\chi = \begin{pmatrix} 0_{(d+1) \times N} \\ 0_{(N-d) \times d} & I_{(N-d) \times (N-d)} \end{pmatrix}; \quad \chi_I = \begin{pmatrix} I_{(d+1) \times N} \\ 0_{(N-d) \times d} & I_{(N-d) \times (N-d)} \end{pmatrix};$$

$$\bar{C}_c = \begin{pmatrix} K_{PD}\left(1 + \dfrac{T_d}{T}\right) & -K_{PD}\dfrac{T_d}{T} \\ K_{PD} & 0 \\ \vdots & \vdots \\ K_{PD} & 0 \end{pmatrix}_{(N+1) \times 2},$$ (4.5b)

And finally, when $\tau_k = dT$, $d = 1..N-1$

$$\chi = \begin{pmatrix} 0_{d \times N} \\ 0_{(N-d) \times d} & I_{(N-d) \times (N-d)} \end{pmatrix}; \quad \chi_I = \begin{pmatrix} I_{d \times N} \\ 0_{(N-d) \times d} & I_{(N-d) \times (N-d)} \end{pmatrix};$$

$$\bar{C}_c = \begin{pmatrix} K_{PD}\left(1 + \dfrac{T_d}{T}\right) & -K_{PD}\dfrac{T_d}{T} \\ K_{PD} & 0 \\ \vdots & \vdots \\ K_{PD} & 0 \end{pmatrix}_{N \times 2},$$ (4.5c)

Note that, as a delay-independent control solution, the controller in (4.3)-(4.5c) is defined irrespective of the delay $\tau_k$, that is, no controller's parameter is retuned according to the delay (as needed, for example, in [7]). In addition, when dropouts occur ($d_k^{lr} = 0$ or $d_k^{rl} = 0$), as the configuration of $C_c$ does not depend on the delay because of the uniform actuation of the estimated control signal, the closed-loop model is also not dependent on the delay, $A_{cl,0}$. However, when no packet dropout occurs ($d_k^{lr} = d_k^{rl} = 1$), the configuration of $\bar{C}_c$ does depend on the time-varying delay in order to satisfy the consequent input pattern to the plant. In this case, the plant is subjected to variations in the instants where the input commands are presented, and hence, matrices $A_P$, $B_P$, $C_P$ in (4.1) may vary from sensor period to sensor period. This leads to represent the NCS closed-loop model as an LTV system depending on $\tau_k$, $A_{cl,1}(\tau_k)$.

To obtain the closed-loop system, the following definition is used

$$w_k^{NT} = y_k^{NT} - \hat{y}_k^{NT} = C_p\left(x_k^{NT} - \hat{x}_k^{NT}\right), \tag{4.6}$$

and this state variable is adopted $\tilde{x}_k^{NT} = \begin{pmatrix} x_k^{NT} & \varphi_k^{NT} & w_k^{NT} \end{pmatrix}^{\mathrm{T}}$. Some manipulations lead to the closed-loop system:

$$\tilde{x}_{k+1}^{NT} = \begin{cases} A_{cl,0}\tilde{x}_k^{NT} & , \text{if } d_k^{lr} = 0 \text{ or } d_k^{rl} = 0 \\ A_{cl,1}(\tau_k)\tilde{x}_k^{NT} & , \text{if } d_k^{lr} = d_k^{rl} = 1 \end{cases}, \tag{4.7}$$

where

$$A_{cl,0} = \begin{bmatrix} A_p & B_p C_c & 0 \\ \bar{B}_c C_p - B_c C_p A_p & A_c - B_c C_p B_p C_c & -\bar{B}_c - A_p B_c \\ 0 & 0 & A_p \end{bmatrix},$$

$$A_{cl,1}(\tau_k) = \begin{bmatrix} A_p & B_p \bar{C}_c & 0 \\ \bar{B}_c C_p - B_c C_p A_p & A_c - B_c C_p B_p \bar{C}_c & 0 \\ 0 & 0 & 0 \end{bmatrix}, \tag{4.8}$$

*4.2 Closed-loop stability*

To assess the closed-loop system stability, the next Theorem can be enunciated.

*Theorem:* Given $P[\tau_k]$ in (2.2), which is normalized to take values in $[0,1)$, the closed-loop system in (4.7)-(4.8) is stable if there exists a positive definitive solution $Q = Q^{\mathrm{T}} > 0$ for the following LMIs

$$\begin{aligned} A_{cl,0}^{\mathrm{T}} \cdot Q \cdot A_{cl,0} - Q &< 0 \\ \sum_{j=1}^{l} P[\vartheta_j] A_{cl,1}^{\mathrm{T}}(\vartheta_j) \cdot Q \cdot A_{cl,1}(\vartheta_j) - Q &< 0 \end{aligned}, \tag{4.9}$$

where $\vartheta$ is a dummy parameter ranging in a set $\Theta$ where $\tau_k$ is assumed to take values in, being $\Theta$ an interval $[0, \tau_{\max}]$, and $l$ the number of equally spaced values to get a dense enough gridding in $\vartheta$. To solve (4.9), widely known methods [3] can be used.

*Proof:* Let $\tilde{V}_k^{NT} = \left(\bar{x}_k^{NT}\right)^{\mathrm{T}} \cdot Q \cdot \bar{x}_k^{NT}$ be a Lyapunov candidate. Taking $E\left[\tilde{V}_k^{NT}\right]$ as the statistical expectation for the Lyapunov function, and assuming a probabilistic LMI gridding procedure, then the expectation of the increment $E\left[\Delta \tilde{V}_k^{NT}\right]$ along subsystem $A_{cl,1}(\tau_k)$ in (4.7)-(4.8) can be obtained as follows

$$\begin{aligned}
E\left[\Delta \tilde{V}_k^{NT}\right] &= E\left[\tilde{V}_{k+1}^{NT} - \tilde{V}_k^{NT}\right] = E\left[\left(\tilde{x}_{k+1}^{NT}\right)^T \cdot Q \cdot \tilde{x}_{k+1}^{NT} - \left(\tilde{x}_k^{NT}\right)^T \cdot Q \cdot \tilde{x}_k^{NT}\right] = \\
&= E\left[\left(A_{cl,1}(\tau_k) \cdot \tilde{x}_k^{NT}\right)^T \cdot Q \cdot A_{cl,1}(\tau_k) \cdot \tilde{x}_k^{NT} - \left(\tilde{x}_k^{NT}\right)^T \cdot Q \cdot \tilde{x}_k^{NT}\right] = \\
&= E\left[\left(\tilde{x}_k^{NT}\right)^T A_{cl,1}^T(\tau_k) \cdot Q \cdot A_{cl,1}(\tau_k) \cdot \tilde{x}_k^{NT} - \left(\tilde{x}_k^{NT}\right)^T \cdot Q \cdot \tilde{x}_k^{NT}\right] = \\
&= E\left[\left(\tilde{x}_k^{NT}\right)^T \left(A_{cl,1}^T(\tau_k) \cdot Q \cdot A_{cl,1}(\tau_k) - Q\right) \cdot \tilde{x}_k^{NT}\right] = \\
&= \left(\tilde{x}_k^{NT}\right)^T \left(\sum_{j=1}^l P\left[\vartheta_j\right] A_{cl,1}^T(\vartheta_j) \cdot Q \cdot A_{cl,1}(\vartheta_j) - Q\right) \tilde{x}_k^{NT} < 0
\end{aligned} \qquad (4.10)$$

By including every possible delay in $\vartheta$, system stability can be ensured for the different configurations of the controller in (4.5a)-(4.5c). Regarding $A_{cl,0}$, a similar development can be carried out, not considering neither the expectation nor the delay. Finally, both developments lead to (4.9).

*Discussion on feasibility:* Note that the first inequality in (4.9), i.e. $A_{cl,0}^T \cdot Q \cdot A_{cl,0} - Q < 0$, would never hold if the plant were open-loop unstable (as the one considered in next Sections). However, if no model uncertainties and an accurate prediction were assumed (i.e. $w_k^{NT} = 0$), the designed controller would be able to stabilize the plant. Then, under these conditions, in order to assess feasibility for the LMIs in (4.9), $A_{cl,0}$ may be replaced by

$$A_{cl,0} = \begin{bmatrix} A_p & B_p C_c \\ \bar{B}_c C_p - B_c C_p A_p & A_c - B_c C_p B_p C_c \end{bmatrix}, \qquad (4.11)$$

And, $A_{cl,1}(\tau_k)$ may consequently reduce its dimensions as follows

$$A_{cl,1}(\tau_k) = \begin{bmatrix} A_p & B_p \bar{C}_c \\ \bar{B}_c C_p - B_c C_p A_p & A_c - B_c C_p B_p \bar{C}_c \end{bmatrix}, \qquad (4.12)$$

# 5. Simulation results

This section is split into three subsections. In subsection 5.1, the data used in the simulations will be presented. In addition, stability for the proposed NCS will be assessed by means of the LMIs stated in Section 4. The aim of subsection 5.2 will be to reveal the benefits of the proposed control solution by comparison with a delay-dependent approach in [7]. The system responses are shown and analyzed via some cost indexes. Finally, in subsection 5.3, model mismatches are considered for the delay-independent control solution, and their consequent effects on the system response are depicted and analyzed by means of some cost indexes.

*5.1. Simulation data. Control system stability assessment*

The process to be controlled is a Cartesian robot manufactured by Inteco, specifically, the 3D CRANE module (see in Figure 4). The rail measures of this plant for each axis are: X=0.050m, Y=0.040m, Z=0.050m.

Focusing on the X axis, its identification leads to

$$G_p(s) = \frac{6.3}{s(s+17.7)} \text{ m/c.a.u.}, \quad (5.1)$$

where c.a.u. means control action units, which are generated by a PWM signal normalized in the range [0,1].

The system also presents two non-linear behaviors to be taken into account in real-time implementation: saturation limits of control actions in ±1, and dead zone values of ±0.06. Both of them are identified experimentally and measured in normalized c.a.u.

For a fair comparison in the next subsection, the simulations are based on the consideration of a UDP network with the same parameters used in [7], that is,

- the network-induced delays are given by the histogram shown in Figure 5, which can be modeled as a generalized exponential distribution (2.2), where $\tau_k$ takes values in the range $\Theta = [0.04, 0.08)$. As $\tau_{max} = 0.08$, then the sensor period can be chosen, for instance, as $NT=0.2$s in order to avoid packet disorder.

- the control design considers a conventional PID controller, being $K_p = 12$, $T_d = 0.01$ and $T_i = 3.5$ in (3.2b) in order to reach the following specifications: null steady-state error, settling time around 4s, and overshot around 5%. From this controller, the dual-rate control is set up using (3.3)-(3.12), where $N=2$ is assumed.

- the packet dropouts are modeled by means of a Bernoulli distribution (2.3), being $p = p^{lr} = p^{rl} = 0.3$ and $M=3$.

In order to assess the stability of the setup in probabilistic time-varying delays, the LMI gridding in (4.9) has been carried out taking $l=20$ grid points so as to compute the closed-loop realization $A_{cl,1}(\vartheta_j)$ for the parameter space $\Theta$. From Figure 5, $P[\vartheta_j]$ is normalized in order to take values in $[0,1)$. As the plant in (5.1) is open-loop unstable, both $A_{cl,0}$ and $A_{cl,1}(\vartheta_j)$ in (4.11) and (4.12) are respectively used to check LMI feasibility. As a result, stability for the proposed NCS can be guaranteed, since the following feasible solution $Q$ exists

$$Q = 10^2 \begin{bmatrix} 0.311893 & 0.056373 & -0.076933 & -0.063061 \\ 0.056373 & 0.661763 & 0.347536 & 0.023331 \\ -0.076933 & 0.347536 & 1.247753 & -0.032175 \\ -0.063061 & 0.023331 & -0.032175 & 0.329738 \end{bmatrix}, \quad (5.2)$$

*5.2. System responses. Cost indexes $J_1$ and $J_2$*

In this subsection, the delay-independent control solution proposed in this work will be compared to the delay-dependent approach based on a gain-scheduling technique presented in [7].

Figure 6 shows this comparison, where filtered step references (in dashed line) are used in order to avoid the saturation of the control signal. Note that the sequence of packet dropouts is represented in the time axis in such a way that each point indicates a packet dropout in the time instant where it is plotted. If the point increases its value in the vertical axis, then consecutive dropouts are occurring in this instant. Circles are used to show the delay-independent control solution, and a thin line is used to depict the delay-dependent approach. In addition, two more responses are included in Figure 6: the nominal, desired response, that is, the output obtained for the dual-rate control when no delay and no dropout are considered (in bold line), and the response for the dual-rate control when time-varying delays and dropouts occur but no prediction stage is included (in dotted line).

The simulation results show the next conclusions: i) as expected, the dual-rate control solution with no prediction stage is negatively affected by delays and dropouts, exhibiting the worst behavior, ii) as the delay-dependent approach includes a prediction stage, it is able to restore the control performance, but it is not able to accurately reach the desired specifications, and iii) the delay-independent control solution is able to achieve the nominal, desired behavior.

In order to better quantify these results, the cost indexes $J_1$ and $J_2$ are stated. $J_1$ is based on the Integral of Absolute Error (IAE), and $J_2$ on the overshoot value.

Let us consider the array $Y$, which includes the four control responses shown in Figure 6, that is, $Y = [Y_{Nom}, Y_{NP}, Y_{DD-P}, Y_{DI-P}]$, being $Y(1) = Y_{Nom}$ the output for the nominal (no-delay, no-dropout) dual-rate control, $Y(2) = Y_{NP}$ the output for the dual-rate control with no prediction stage and occurring delays and dropouts, $Y(3) = Y_{DD-P}$ the output for the delay-dependent approach (which includes prediction stage), and $Y(4) = Y_{DI-P}$ the output for the delay-independent proposal (which also includes prediction stage). From $Y$, the following accumulated (integrated) error $E_Y$ in a range of time instants $\Gamma$ can be computed

$$E_Y(i) = \sum_{\Gamma} |Y(i) - Y_{Nom}|, \quad i = 1..4, \tag{5.3}$$

Then, the $J_1$ cost index takes this form

$$J_1(i) = 100 - \frac{E_Y(i)}{E_Y(2)} 100 \, (\%), \quad i = 1..4, \tag{5.4}$$

being $E_Y(2)$ the worst expected accumulated error, that is, the error calculated for $Y(2) = Y_{NP}$. Therefore, the rest of the errors are measured by $J_1$ as an improvement (in %) with respect to $E_Y(2)$.

Additionally, from $Y$, the following overshoot $O_Y$ in a range of time instants $\Gamma$ can be calculated (considering positive -max- or negative -min- filtered step references)

$$O_Y(i) = \max\left(\left|\max_\Gamma Y(i) - \max_\Gamma Y_{Nom}\right|, \left|\min_\Gamma Y(i) - \min_\Gamma Y_{Nom}\right|\right), \quad i=1..4, \tag{5.5}$$

Then, the $J_2$ cost index is defined as

$$J_2(i) = 100 - \frac{O_Y(i)}{O_Y(2)} 100 \ (\%), \quad i=1..4, \tag{5.6}$$

being $O_Y(2)$ the worst expected overshoot, that is, the overshoot obtained for $Y(2) = Y_{NP}$. Similarly to $J_1$, the rest of the overshoots are measured by $J_2$ as an improvement (in %) with regard to $O_Y(2)$.

Table 1 summarizes the cost indexes $J_1$ and $J_2$ obtained for each output. The delay-dependent control solution significantly improves $J_2$ (around 85%) with respect to the worst response $Y_{NP}$, but it exhibits a poor improvement (around 42%) for the value $J_1$. Nevertheless, the delay-independent control approach is able to accurately achieve the same control properties as the nominal dual-rate control solution, since $J_1$ and $J_2$ practically reach 100%.

## 5.3. Model mismatch. Cost indexes $J_3$ and $J_4$.

As the proposed delay-independent control solution is model-based, both the controller design and the prediction computation depend on how accurate the model represents the plant behavior. If some uncertainty between plant and model were assumed, the robustness of the model-based control proposal could be checked.

In this subsection, ~~some~~ a certain model mismatch in the characteristic parameters of the plant, say, the static gain $K$ and the time constant $\tau$, will be considered. Let us consider a percentage $q\%$ of increase in $K$, $q\% \cdot \Delta K$, and a percentage $r\%$ of increase in $\tau$, $r\% \cdot \Delta\tau$. Figure 7 shows a comparison among the nominal response and the outputs obtained when $[20\% \cdot \Delta K, \ 8\% \cdot \Delta\tau]$ and [$30\% \cdot \Delta K$, $12\% \cdot \Delta \tau$]. As expected, the larger the percentage of uncertainty is considered, the worse the behavior becomes (with regard to the nominal response). Despite assuming a significant model mismatch when [$30\% \cdot \Delta K$, $12\% \cdot \Delta \tau$], the system remains stable. However, the response is clearly worsened, exhibiting a 60% increase in settling time, and a 6% increase in overshoot.

In order to better quantify these results and the robustness of the approach, the cost indexes $J_3$ and $J_4$ are stated. $J_3$ is based on the Integral of Absolute Error (IAE), and $J_4$ on the overshoot value.

Let us consider the matrix $W$, which includes different outputs for the proposed delay-independent approach. These outputs are obtained as a result of varying $q\%\Delta K$ and $r\%\Delta \tau$. In this study, $q$ takes the values $q$=0, 20, 30, and $r$ the values $r$=0, 8, 12. Combining these values, nine different responses in $W$ can be considered. The nominal response $Y_{Nom}$ is obtained when $q$=$r$=0. As previously commented, the worst behavior will be obtained for $q$=30 and $r$=12, since it represents the largest model mismatch. Let us assume this behavior as the worst permissible one. From $W$ the following accumulated (integrated) error $E_W$ in a range of time instants $\Gamma$ can be computed

$$E_W(i_r,i_q) = \sum_{\Gamma} \left| W(i_r,i_q) - Y_{Nom} \right|, \quad i_r, i_q = 1..3, \tag{5.7}$$

Then, the $J_3$ cost index takes this form

$$J_3(i_r,i_q) = 100 - \frac{E_W(i_r,i_q)}{E_W(3,3)} 100 \, (\%), \quad i_r, i_q = 1..3, \tag{5.8}$$

being $E_W(3,3)$ the worst permissible accumulated error, that is, the error reached when $r$=12 and $q$=30. Therefore, the rest of the errors are measured by $J_3$ as an improvement (in %) with respect to $E_W(3,3)$.

Additionally, from $W$, the following overshoot $O_W$ in a range of time instants $\Gamma$ can be calculated (considering positive -max- or negative -min- filtered step references)

$$O_W(i_r,i_q) = \max\left(\left|\max_\Gamma W(i_r,i_q) - \max_\Gamma Y_{Nom}\right|, \left|\min_\Gamma W(i_r,i_q) - \min_\Gamma Y_{Nom}\right|\right), \quad i_r, i_q = 1..3, \qquad (5.9)$$

And then, the $J_4$ cost index is defined as

$$J_4(i_r,i_q) = 100 - \frac{O_W(i_r,i_q)}{O_w(3,3)} 100 \, (\%), \quad i_r, i_q = 1..3, \qquad (5.10)$$

being $O_W(3,3)$ the worst permissible overshoot, that is, the overshoot reached when $r=12$ and $q=30$. Similarly to $J_3$, the rest of the overshoots are measured by $J_4$ as an improvement (in %) with regard to $O_W(3,3)$.

Tables 2 and 3 respectively summarize the cost indexes $J_3$ and $J_4$ obtained for each output depending on the model mismatch considered. As expected, the smaller the percentage of mismatch is considered, the larger $J_3$ and $J_4$ become, that is, a closer behavior to the nominal one is obtained.

# 6. - Experimental results

To validate the simulation results obtained in section 5, a laboratory test-bed platform is set up, which includes the CRANE module previously presented, two computers and an Ethernet cable.

One computer is directly connected to the plant and includes the local part of the control system. The aims of this computer are: firstly, to be in charge of the sampling measurement and transmission at $NT=0.2$s; secondly, to receive the packet which includes the current and predicted PI control signals; thirdly, to compute, and inject to the plant, the fast-rate PD control actions at $T=0.1$s.

The second computer performs the remote part of the controller, receiving the outputs of the plant, calculating the current and the predicted slow-rate PI control actions, and sending back these actions to the local system. Both computers are connected by a UDP network through an Ethernet cable that performs the local-to-remote and remote-to-local links. In order to obtain the same conditions as those considered in simulation, packet delays and packet dropouts are modified by software.

Figure 8 shows the outputs obtained in the experiment, which clearly validates the trend observed in Figure 6, that is, the delay-independent control solution improves the results obtained by the delay-dependent approach in [7], being able to reach the nominal, desired behavior. To better validate the results, Table 4 details the cost indexes $J_1$ and $J_2$ computed for the experiment. While the delay-independent approach exhibits values for $J_1$ and $J_2$ very close to 100%, which means that it practically achieves the nominal behavior, the delay-dependent proposal presents percentages which are significantly smaller.

Finally, to strengthen the previous conclusions, Figure 9 compares the behavior achieved by every control solution when the robot tries to track a 2D trajectory based on the well-known Lissajous curves (see, e.g., [1]). As expected, the response for the dual-rate control system when no prediction stage is included presents the worst behavior, mainly, when the curves are more pronounced. Once again, the difference between the delay-dependent strategy in [7] and the delay-independent control solution proposed in this work is clearly observed, since, whereas the former does not accurately track the nominal response, the latter does.

## 7. - Conclusions

In this work, an NCS is presented where time-varying delays, packet dropouts and packet disorder can occur. By means of a novel delay-independent control solution, which integrates packet-based control strategies with dual-rate and predictor-based control techniques, the above problems are faced, and in addition, ~~the~~ usage of the network resources is reduced, while keeping the desired control specifications.

Control system stability is ensured in terms of LMIs. The benefits of the control solution are illustrated by simulation, and experimentally validated by means of a robotic test-bed platform.

# Acknowledgments

This work is funded by European Commission as part of Project H2020-SEC-2016-2017, Topic: SEC-20-BES-2016 – Id: 740736 – "C2 Advanced Multi-domain Environment and Live Observation Technologies" (CAMELOT). Part WP5 supported by Tekever ASDS,Thales Research & Technology, Viasat Antenna Systems, Universitat Politècnica de València, Fundação da Faculdade de Ciências da Universidade de Lisboa, Ministério da DefesaNacional – Marinha Portuguesa, Ministério da Administração Interna Guarda Nacional Republicana.

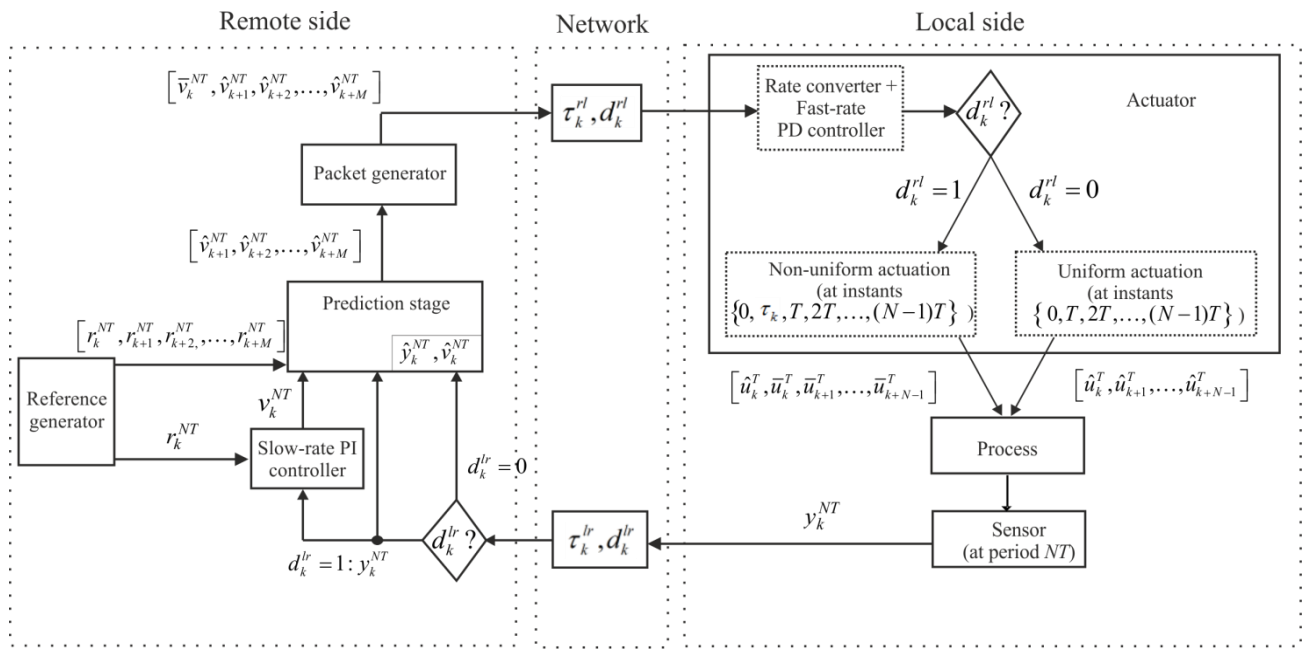

Figure 1. NCS scenario.

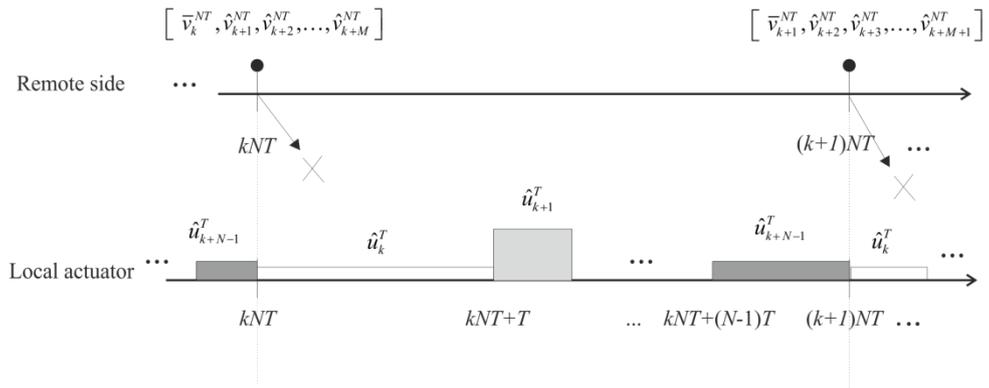

Figure 2. Packet dropout ($d_k^{rl} = 0$).

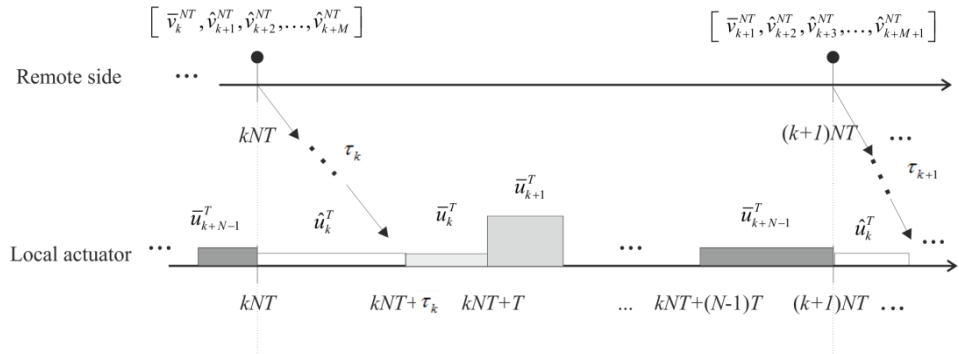

Figure 3. No packet dropout ($d_k^{rl} = 1$).

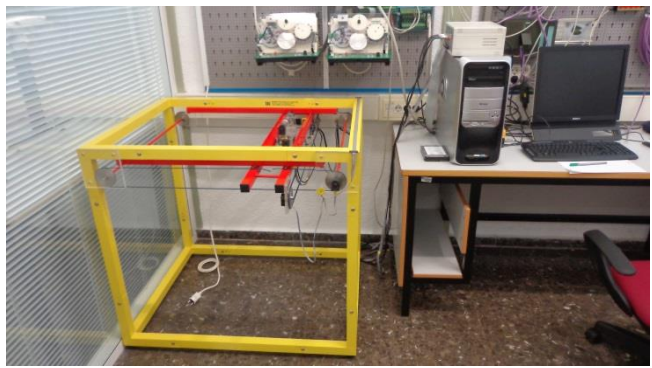

Figure 4. Cartesian robot (3D CRANE module).

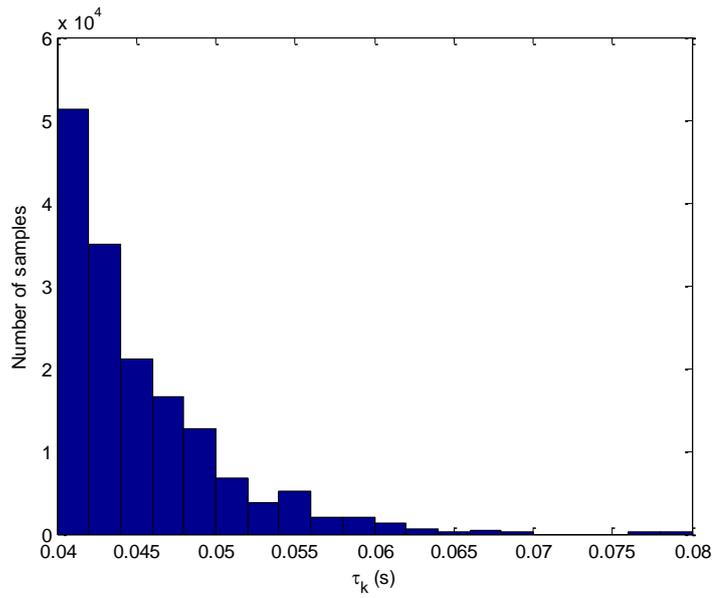

Figure 5. Delay histogram.

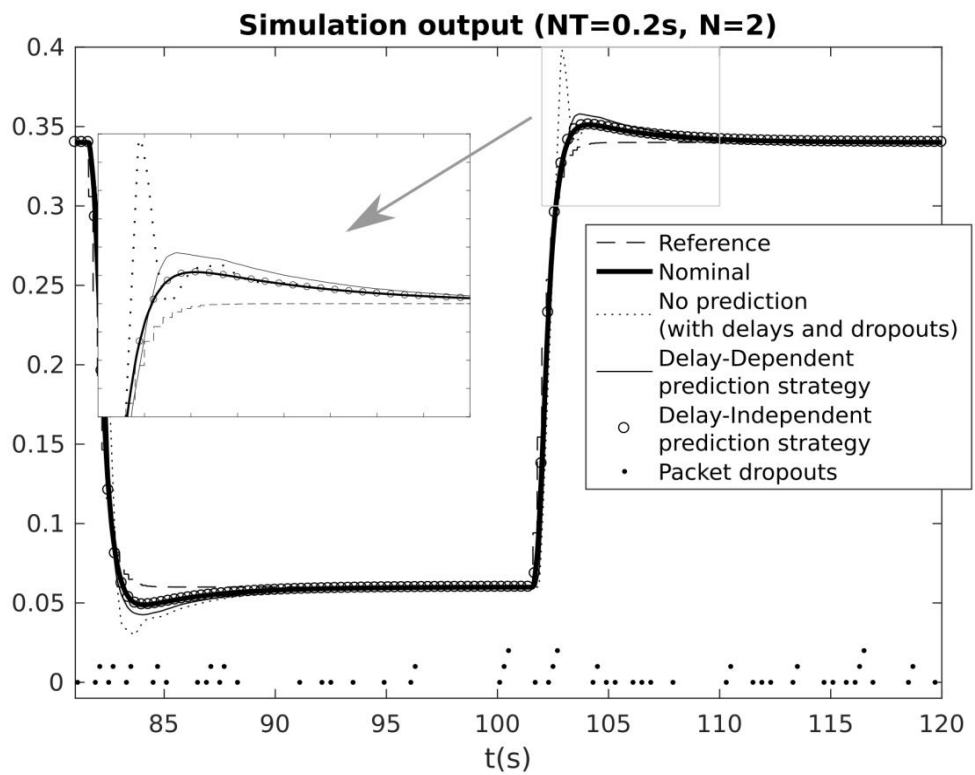

Figure 6. Comparison: nominal vs no prediction vs delay-dependent vs delay-independent.

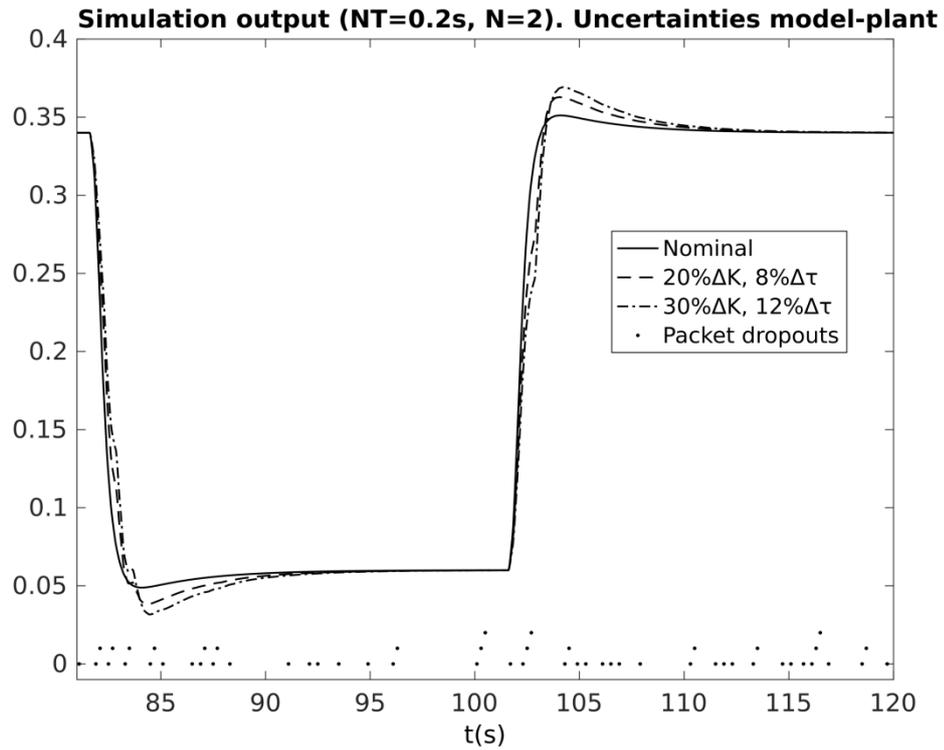

Figure 7. Comparison: nominal vs delay-independent (model mismatch).

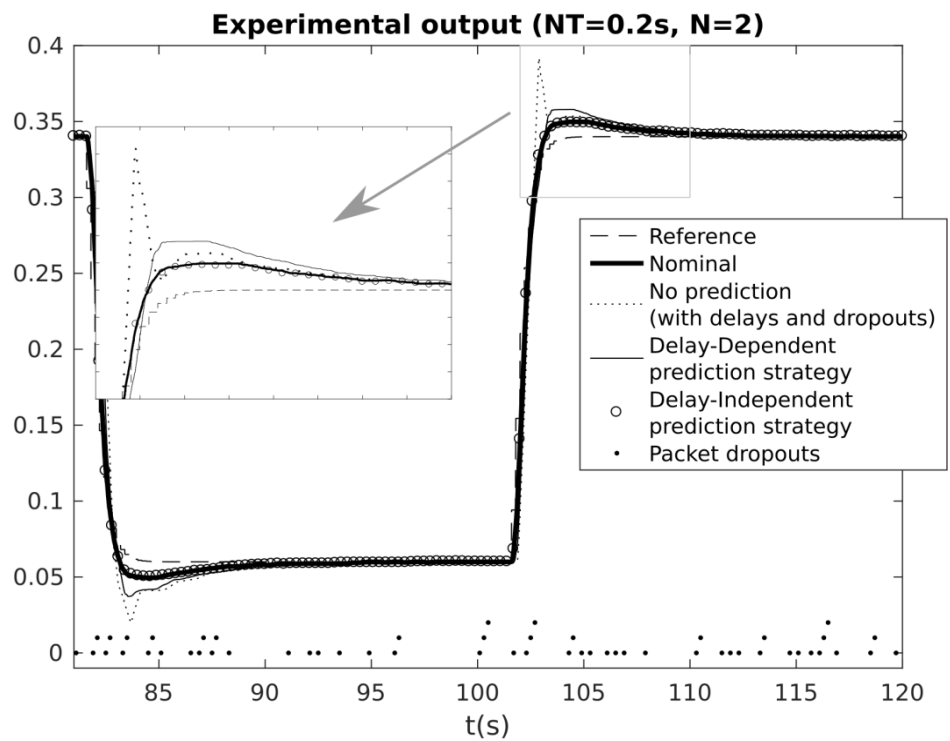

Figure 8. Comparison: nominal vs no prediction vs delay-dependent vs delay-independent.

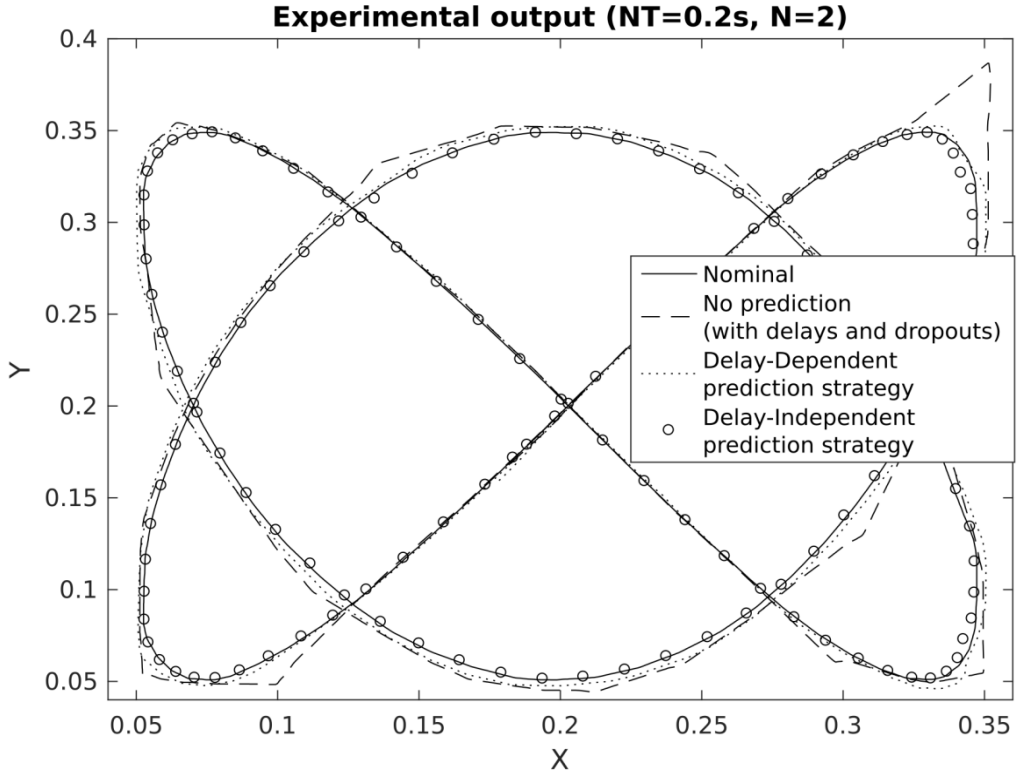

Figure 9. Comparison (Lissajous curves): nominal vs no prediction vs delay-dependent vs delay-independent.

| Output | $E_Y$ | $J_1(\%)$ | $O_Y$ | $J_2(\%)$ |
|---|---|---|---|---|
| $Y_{NP}$ | 291.19 | 0 | 0.047 | 0 |
| $Y_{DD-P}$ | 168.64 | 42.09 | 0.007 | 85.64 |
| $Y_{DI-P}$ | 0.02 | 99 | 0 | 100 |

Table 1. Simulation: accumulated error $E_Y$ and cost index $J_1$; overshoot $O_Y$ and cost index $J_2$.

| $E_W$ | | $q\%\Delta K$ | | |
|---|---|---|---|---|
| | | 0 | 20 | 30 |
| | 0 | 0 | 183.46 | 302.60 |
| $r\%\Delta\tau$ | 8 | 59.82 | 248.62 | 373.11 |
| | 12 | 88.56 | 280.65 | 407.19 |

| $J_3$ | | $q\%\Delta K$ | | |
|---|---|---|---|---|
| | | 0 | 20 | 30 |
| | 0 | 100 | 54.94 | 25.68 |
| $r\%\Delta\tau$ | 8 | 85.31 | 38.94 | 8.36 |
| | 12 | 78.25 | 31.07 | 0 |

Table 2. Accumulated error $E_W$ and cost index $J_3$.

| $O_W$ | | $q\%\Delta K$ | | |
|---|---|---|---|---|
| | | 0 | 20 | 30 |
| | 0 | 0 | 0.0087 | 0.0142 |
| $r\%\Delta\tau$ | 8 | 0.0027 | 0.0117 | 0.0169 |
| | 12 | 0.0041 | 0.0130 | 0.0182 |

| $J_4$ | | $q\%\Delta K$ | | |
|---|---|---|---|---|
| | | 0 | 20 | 30 |
| | 0 | 100 | 52.19 | 21.97 |
| $r\%\Delta\tau$ | 8 | 85.16 | 35.71 | 7.14 |
| | 12 | 77.47 | 28.57 | 0 |

Table 3. Overshoot $O_W$ and cost index $J_4$.

| **Output** | $E_Y$ | $J_1(\%)$ | $O_Y$ | $J_2(\%)$ |
|---|---|---|---|---|
| $Y_{NP}$ | 290.55 | 0 | 0.042 | 0 |
| $Y_{DD-P}$ | 207.53 | 28.58 | 0.012 | 70.80 |
| $Y_{DI-P}$ | 19.53 | 93.28 | 0.001 | 97.62 |

Table 4. Experiment: accumulated error $E_Y$ and cost index $J_1$; overshoot $O_Y$ and cost index $J_2$.